\documentclass[11pt,a4paper]{article}
\usepackage[utf8]{inputenc}
\usepackage{amsmath}
\usepackage{amsfonts}
\usepackage{amssymb}
\usepackage{bm}
\usepackage{pdflscape}
\usepackage[margin=1in,footskip=0.25in]{geometry}
\usepackage{graphicx}
\usepackage{multirow}
\usepackage{graphicx}
\usepackage{footnote}
\usepackage{lscape}
\usepackage[table,xcdraw]{xcolor}
\usepackage{hyperref}
\usepackage{float}

\usepackage{longtable}
\usepackage{courier}
\usepackage{hyperref}
\usepackage[printwatermark]{xwatermark}
\begin{document}
\noindent\textbf{\Large{Estimating the proportion of true null hypotheses with application in microarray data}}
\\
\begin{center}
    \textbf{(Under revision)}
\end{center}
\begin{center}
Aniket Biswas\\
Department of Statistics\\
Dibrugarh University\\
Dibrugarh, Assam, India-786004\\
Email:\textit{biswasaniket44@gmail.com}
\end{center}

%\doublespacing
\section*{Abstract}
A new formulation for the proportion of true null hypotheses $(\pi_0)$, based on the sum of all $p$-values and the average of expected $p$-value under the false null hypotheses has been proposed in the current work. This formulation of the parameter of interest $\pi_0$ has also been used to construct a new estimator for the same. The proposed estimator removes the problem of choosing tuning parameters in the existing estimators. Though the formulation is quite general, computation of the new estimator demands use of an initial estimate of $\pi_0$. The issue of choosing an  appropriate initial estimator is also discussed in this work. The current work assumes normality of each gene expression level and also assumes similar tests for all the hypotheses. Extensive simulation study shows that, the proposed estimator performs better than its closest competitor, the estimator proposed in Cheng et al., 2015 over a substantial continuous subinterval of the parameter space, under independence and weak dependence among the gene expression levels. The proposed method of estimation is applied to two real gene expression level data-sets and the results are in line with what is obtained by the competing method.\\ 

\noindent \textbf{Keywords}: True null, p-value, Normality, t-test, Effect size, Microarray data, Expected p-value.\\

\noindent \textbf{MS 2010 classification}: 62F10, 62P10. 

\section{Introduction}
In this era of high throughput devices, huge datasets are easily available to answer complicated decision-making questions. In microarray experiments, data on thousands of genes are available and from that large number of genes, the task is to identify the differentially expressed genes between a set of control subjects and a set of treatment subjects for making further scientific experimentation efficient. Thus, testing thousands of hypotheses simultaneously and being able to make higher number of rejections with control over false discovery rate (FDR) (see Benjamini and Hochberg, 1995) is desirable. Benjamini-Hochberg procedure for controlling the FDR is originally conservative. A reliable estimate of $\pi_0$ can be used to eliminate conservative bias of the same (see Benjamini and Yekutieli, 2001). Efficient estimation of $\pi_0$ can improve algorithms, controlling family wise error rate through reduction in false negative rate (see Hochberg and Benjamini, 1990; Finner and Gontscharuk, 2009).\\

For empirical Bayesian motivation behind FDR given in Storey, 2002, $m$ related but independent hypotheses viz. $H_1, H_2, ..., H_m$ are considered. The $i$-th null hypotheses $H_i$ is seen as an indicator variable. Here, $H_i=1$ indicates that the $i$-th null hypothesis is true and $H_i=0$ indicates the same to be false, for all $i\in \mathcal{I}=\{1, 2, ..., m\}$. It is to be realized that the hypotheses are no longer considered only as partitions of the parameter space. In fact, $H_i$'s are Bernoulli random variables with success probability $\pi_0\in (0,1)$. Hence, $m_0=\sum_{i=1}^m H_i$, number of true null hypotheses is a binomial random variable with index $m$ and parameter $\pi_0$. Unfortunately, $H_i$'s and thus $m_0$ remain unrealized in a given multiple testing scenario. Usually $p$-values are considered as test-statistics since one gets similar critical region for each test, whatever be the nature of hypothesis to be tested. Throughout this article, $p$-values are denoted by $p$ irrespective of it being an observed value or a random variable. The notation holds the meaning in accordance with the situation. From this motivation, Langaas et al.(2005) put forward the following two-component mixture model for marginal density of $p$-value 
\begin{equation}\label{eq:twocomptpdens}
    f(p)=\pi_0\,f_0(p)\,+\, (1-\pi_0)\,f_1(p)\quad \textrm{for}\quad 0<p<1.
\end{equation}
Here, $f_0$ and $f_1$ denote $p$-value densities under null and alternative hypothesis, respectively. The current work assumes that the null hypotheses to be tested are simple and the corresponding test-statistics are absolutely continuous, which is quite common in the existing literature on estimating $\pi_0$. Thus, $p$ is distributed uniformly over $(0,1)$ under the null hypothesis. Under alternative hypothesis, $p$ is stochastically smaller than the uniform variate and $f_1(p)\to 0$ for $p$ approaching unity. Let $\mathcal{T}$ denote $\{i\in \mathcal{I}: H_i=1\}$, the set of indices corresponding to the originally true null hypotheses. Obviously, cardinality of $\mathcal{T}$ is $m_0$. Also let $\mathcal{F}=\mathcal{I}-\mathcal{T}$ and thus, cardinality of $\mathcal{F}$ is $m_1=m-m_0$. The set of $p$-values, $\{p_1, p_2, ..., p_m\}$ is obtained by performing appropriate test for each of the hypotheses where, $p_i$ denotes the $p$-value corresponding to $H_i$ for all $i\in \mathcal{I}$.\\

For Storey's estimator (see Storey, 2002), existence of a tuning parameter $\lambda\in(0,1)$ is assumed such that, $f_1(p)=0$ for $p\geq \lambda$. For such fixed choice of $\lambda$, $W(\lambda)=\sum_{i\in\mathcal{I}}I(p_i\geq\lambda)=W_1(\lambda)+W_0(\lambda)$, where $W_1(\lambda)=\sum_{i\in \mathcal{T}}I(p_i\geq\lambda)$ and $W_0(\lambda)=\sum_{i\in \mathcal{F}}I(p_i\geq\lambda)$ ($I$ denotes the indicator function). Since $f_1(p)=0$ for $p\geq \lambda$, $W_0(\lambda)=0$ and hence, $E(W(\lambda))=m\pi_0(1-\lambda)$. Thus for a fixed choice of $\lambda$, an estimator for $\pi_0$ is $\hat{\pi}_0(\lambda)=W(\lambda)/[m(1-\lambda)]$. For a subjectively chosen $\Lambda=\{0, 0.05, 0.10, ..., 0.95\}$, Storey(2002) and Storey et al.(2004) discuss a bootstrap routine to find the best choice of $\lambda\in \Lambda$, viz. $\lambda_{best}$. Thus, Storey's bootstrap estimator is $\hat{\pi}_0^B=\hat{\pi}_0(\lambda_{best})$. This estimator has an inherent upward bias due to the crucial assumption $W_0(\lambda)=0$. Cheng et al.(2015) worked out the bias in $\hat{\pi}_0^B$ and formulated $\pi_0= [E(W(\lambda))-mQ(\lambda)]/[m(1-\lambda)-mQ(\lambda)]$. Here, $Q(\lambda)=\sum_{i\in \mathcal{F}}Q_{\delta_i}(\lambda)$ is the average upper tail probability of $p$-values corresponding to $H_i$'s for $i\in \mathcal{F}$ and $\delta_i$'s are effect sizes of the same. $\hat{Q}(\lambda)$ being an estimator for $Q(\lambda)$, the plug-in estimator for $\pi_0$ from Cheng's formulation is $\hat{\pi}_0^U(\lambda)=[W(\lambda)-m\hat{Q}(\lambda)]/[m(1-\lambda)-m\hat{Q}(\lambda)]$. For obtaining $\hat{Q}(\lambda)$, an initial estimate of $\pi_0$ is needed. Cheng et al.(2015) suggested to use $\hat{\pi}_0^B$ as initial estimator for computing $\hat{\pi}_0^U(\lambda)$. For reduction of variance, $\Lambda=\{0.20, 0.25, ..., 0.5\}$ is considered in the same spirit as in Jinag and Doerge(2008). The final bias and variance reduced estimator for $\pi_0$ is 
\begin{equation}
    \hat{\pi}_0^U=\frac{1}{\# \Lambda}\sum_{\lambda_j\in \Lambda} min\{1,max\{0,\hat{\pi}_0^U(\lambda_j)\}\}
\end{equation}
where, $\# \Lambda$ denotes the cardinality of $\Lambda$. This approach reduces the  conservative bias in $\hat{\pi}_0^B$. It is worth mentioning that, obtaining $\hat{Q}(\lambda)$ and hence $\hat{\pi}_0^U$ is possible under some practical assumptions. For all the hypotheses to be tested, common tests are to be performed. For example, two sample two-sided t-tests are usually performed for all the genes of microarray gene expression datasets to identify differentially expressed genes. The observations used for each test are assumed to follow a known parametric distribution. For microarray datasets, normality of gene expression levels is a common assumption. The corresponding test-statistics should have a known exact distribution. This assumption enables one to calculate the exact $p$-value for each test. Application of $t$-test takes care of this assumption. The distribution of $p$-values under alternative are labelled by unknown effect sizes $\delta_i$ for $i\in \mathcal{F}$. This issue is discussed in detail in section 3, for three different testing scenarios, viz. single sample Z-test, single sample two-sided t-test and two sample two-sided t-test. As the assumptions are quite practical in nature, model based bias correction methods may turn out to be really efficient. Often such bias corrected estimators suffer from over-correction as pointed out in Cheng et al., 2015 for the estimator proposed in Qu et al., 2012.\\

Note that, $\hat{\pi}_0^U$ depends on the subjectively chosen index set $\Lambda$ and only utilizes $p$-values greater than a specific threshold. In the current work, an estimator is proposed which does not require such subjective choice of tuning parameters. This is achieved by going a step ahead of the work in Cheng et al., 2015 . Expected $p$-value under the alternative hypothesis is computed using the upper-tail probability of the non-null $p$-value or by directly using the density function of the same. This idea of computing expectation of $p$-value under alternative hypothesis has also been discussed in Hung et al.(1997). The new estimation procedure is proposed under same assumptions as discussed in the last paragraph. The estimator also requires an initial estimate of $\pi_0$. Several variants of the proposed estimator  based on the initial choice of $\pi_0$ are studied in this article. There are several other works on estimation of $\pi_0$, not directly related to the current work. The interested readers are referred to Storey and Tibshirani, 2003, Wang et al., 2011 and the references therein for theoretical developments on the topic. For different applications of the quantity $\pi_0$, one may see Miller et al., 2001 and Turkheimer et al., 2001. \\

In the following section, a new formulation of $\pi_0$ is provided along with the estimation procedure. In section 3, properties of non-null $p$-values are discussed. Performance of variants of the new estimator compared to its closest competitor $\hat{\pi}_0^E$ using simulated datasets is given in section 4. In section 5, the new estimation method is applied to two real life microarray gene expression datasets for obtaining the proportion of non-differentially expressed genes. Limitations of the current study and discussions regarding future scopes are contained in the conclusion of this article.

\section{Proposed method of estimation}
The improvement in $\hat{\pi_0}^U$ for estimating $\pi_0$, over $\hat{\pi_0}^B$ is attributable to appropriate model assumption and use of similar tests for each of  the hypotheses. These assumptions make way for the possible reduction in bias of $\hat{\pi_0}^B$. If we agree upon loosing some generality and in turn gain efficiency, model based formulation of $\pi_0$ has been established to be a simple and effective way out. As mentioned in section 1, $\hat{\pi}_0^U$ utilizes only those $p$-values which are greater than some subjectively chosen $\lambda$. Storey introduced this technique to formulate a robust and conservative estimator for $\pi_0$. Since we are ready to lose some generality, it may be a good idea to use all the available $p$-values for formulating an estimator for $\pi_0$. Obviously, likelihood based approach is an option but unfortunately, this is feasible only when the $p$-values are assumed to be independent. Thus, a moment based estimator for $\pi_0$, with bias correction using appropriate model assumption, is proposed and investigated in this work. As the proposed estimator is based on sum (equivalently, mean) of all the $p$-values, assumption of independence among the $p$-values is not considered necessary. As seen later, performance of the proposed estimator deteriorates under dependence among the $p$-values. Note that,
\begin{equation*}
    E\left(\sum_{i\in \mathcal{I}} p_i\right)= E\left(\sum_{i\in \mathcal{T}} p_i\right)+ E\left(\sum_{i\in \mathcal{F}} p_i\right).
\end{equation*}
Now, $E(p_i|i\in \mathcal{T})=1/2$ and let $e_i$ be $E(p_i|i\in \mathcal{T})$, expected $p$-value when $H_i=0$. Thus, $E\left(\sum_{i\in \mathcal{I}} p_i\right)= m_0/2 + \sum_{i\in \mathcal{F}} e_i$. Denote $1/m_1\sum_{i\in \mathcal{F}}e_i$, average of the expected $p$-values under the alternative hypotheses by $e$. Then, 
\begin{equation*}
     E\left(\sum_{i\in \mathcal{I}} p_i\right)= m_0\left(\frac{1}{2}-e\right)\,+\,m\,e\\
    \implies  E\left(\sum_{i\in \mathcal{I}} p_i\right)\,-\,e=m_0\left(\frac{1}{2}-e\right).
\end{equation*}
Let $\bar{p}$ denote $1/m \sum_{i\in \mathcal{I}} p_i$, mean of the available $p$-values. Thus, from the empirical Bayesian interpretation of $\pi_0$, it can be formulated as
\begin{equation}\label{eq:pi0form}
    \pi_0=\frac{E(p)-e}{0.5-e}.
\end{equation}
As mentioned earlier, $p$-value under the alternative hypothesis is stochastically smaller than $p$-value under the null hypothesis. Hence from \ref{eq:twocomptpdens}, $e<E(p)<0.5$ for $\pi_0\in (0,1)$. Thus, the formulation is reasonable in a sense that, RHS in \ref{eq:pi0form} is in $(0,1)$, the parameter space of $\pi_0$. To construct an estimator for $\pi_0$ from \ref{eq:pi0form}, the unknown quantities $E(p)$ and $e$ are to be estimated. A consistent estimator for $E(p)$ is $\bar{p}$. Let $\Tilde{e}$ be a consistent estimator for $e$. Then, a consistent estimator for $\pi_0$ is $\tilde{\pi}_0=(\bar{p}-\tilde{e})/(0.5-\tilde{e})$. Unfortunately, $\tilde{e}$ can only be conceptualized but not realized in practice. In the next subsection, formulation of $\tilde{e}$ is discussed and a working dummy for $\tilde{e}$ is also proposed.

\subsection{Estimating the average of expected $p$-values under alternative}
As the distribution of $p_i$ is labelled by $\delta_i$ for $i\in \mathcal{F}$, $e_i$ is a function of $\delta_i$, i.e. $e_i=e_{\delta_i}$. Structure of $\delta_i$ depends on the test performed. The explicit expression of effect size for different testing scenarios and their consistent estimation is discussed in section 3. For the present topic under consideration, it is assumed that $\hat{\delta}_i$'s consistently estimate $\delta_i$ for $i\in \mathcal{I}$. Hence, $e_{\hat{\delta}_i}$ is consistent for $e_{\delta_i}$, $i\in \mathcal{I}$. Thus,
\begin{equation}\label{eq:etilde to e}
    \tilde{e}=\frac{1}{m_1}\sum_{i\in \mathcal{F}}\hat{e}_i=\frac{1}{m_1}\sum_{i\in \mathcal{F}}e_{\hat{\delta}_i}\to\frac{1}{m_1}\sum_{i\in \mathcal{F}}e_{\delta_i}= \frac{1}{m_1}\sum_{i\in \mathcal{F}}e_i=e,
\end{equation}
in probability. One can expect that $\tilde{e}$ is close to $e$, at least asymptotically. Different settings for generating simulated data are given in section 4. Finite sample comparison between $\tilde{e}$ and $e$ is done using the simulated data. The results are shown in Figure 1 and the plots justify the claim made above even for small to moderate sample size. This comparative study is possible since $\mathcal{F}$ is known for simulated data. But in practice, $\mathcal{F}$ is unknown and thus, $\tilde{e}$ remains unobserved.\\

Suppose $\hat{\pi}_0^I$, an initial estimate of $\pi_0$ is available. The following algorithm is proposed to obtain $\hat{e}$ which will be used in stead of $\tilde{e}$, in the working estimator.
\begin{center}
    \textbf{Algorithm}\\
    (For obtaining $\hat{e}$)
\end{center}

\begin{itemize}
    \item For a given multiple testing scenario, obtain $\{\hat{\delta}_i:i\in \mathcal{I}\}$.
    \item Using the explicit expression for $e_i=e_{\delta_i}$, obtain $\{\hat{e}_i:i\in \mathcal{I}\}$.
    \item Arrange $\hat{e}_i$'s in increasing order. Denote the $i$-th largest value by $\hat{e}_{(i)}$ for $i\in \mathcal{I}$.
    \item Calculate $d=[m\times (1-\hat{\pi}_0^I)]$ ($[x]$ denotes the largest integer contained in $x$). 
    \item For $\mathcal{D}=\{1, 2, ..., d\}$, calculate $\hat{e}=(1/d)\sum_{i\in \mathcal{D}} \hat{e}_{(i)}$.
\end{itemize}
The role of $\hat{\pi}_0^I$ in obtaining $\hat{e}$ is important. Assume that $\hat{\pi}_0^I\geq \pi_0$, almost surely or equivalently $m_1\geq d$. Clearly 
\begin{equation}\label{eq:ehat less etilde}
 \hat{e}=\frac{1}{d}\sum_{i \in \mathcal{D}}\hat{e}_{(i)}\leq \frac{1}{m_1}\sum_{i=1}^{m_1} \hat{e}_{(i)}\leq \frac{1}{m_1}\sum_{i\in \mathcal{F}} \hat{e}_{i}=\tilde{e}, 
\end{equation}
almost surely. Using \ref{eq:etilde to e} and \ref{eq:ehat less etilde}, $\hat{e}\leq e$, in probability. See Figure 1 for the comparative study on $e$, $\tilde{e}$ and $\hat{e}$ (using different initial estimates) with simulated data. The results justify the theoretical claims made for $\hat{e}$. It is worth mentioning that, we refrain from committing over-correction in bias by using $\hat{e}$. This fact is elaborated in the following subsection.
\subsection{Conservative bias and choice of the initial estimate}
Replacing $\tilde{e}$ by $\hat{e}$ in $\tilde{\pi}_0$, obtain the working estimator $\hat{\pi}_0=(\bar{p}-\hat{e})/(0.5-\hat{e})$. For $\bar{p}<0.5$, 
\begin{equation*}
    \frac{\bar{p}-x}{0.5-x}=1-\frac{0.5-\bar{p}}{0.5-x} 
\end{equation*}
increases with decreasing $x$. Thus, $\hat{\pi}_0\geq \tilde{\pi}_0$. As $\tilde{\pi}_0$ is consistent for $\pi_0$, $\hat{\pi}_0$ dominates true $\pi_0$, asymptotically. This conclusion is desirable as its subsequent use in estimation of FDR, overestimation is preferred to underestimation. The simulation results  presented in section 4 justifies this claim for the proposed estimator (using different initial estimates). As $\bar{p}$ is consistent for $E(p)$, the assumption on $\bar{p}$ being less than $0.5$, mean of uniform$(0,1)$-variate is usually satisfied in practice. Still $\hat{\pi}_0$ is not bound to lie in the interval $[0,1]$. Thus, necessary modification is made to $\hat{\pi}_0$ to get the final estimator as
\begin{equation}\label{eq:finalest}
    \hat{\pi}_0^E=min\left\{1,\,max\left\{\frac{\bar{p}-\hat{e}}{0.5-\hat{e}}\right\}\right\}.
\end{equation}
\indent The conservative bias in $\hat{\pi}_0^E$ depends on the availability of $\hat{\pi}_0^I$, which dominates $\pi_0$ almost surely. Unfortunately, an initial estimator with such a strong theoretical property is not currently available. However, estimators for $\pi_0$ with conservative bias can be used as initial estimator. Two such popular choices are Storey's bootstrap estimator $\hat{\pi}_0^B$ and $\hat{\pi}_0^L$, the convest density estimation based estimator proposed in Langaas et al., 2005. While conservative estimator is preferable due to obvious reason, the candidates for $\pi_0^I$ should not overestimate beyond a reasonable limit. If it does so, then the number of indices in $\mathcal{T}\cap \mathcal{D}$ increases resulting in very small value of $\hat{e}$ compared to $e$. This fact makes the amount of model based bias correction negligible, which is not desirable. Figure 1 clearly shows that, the use of $\hat{\pi}_0^B$ and $\hat{\pi}_0^L$ as $\hat{\pi}_0^I$ for computing $\hat{e}$ results in substantial approximation of $e$. For further understanding of the influence of the initial estimate, one-step iteration of $\hat{\pi}_0^E$ may be performed. That is, in first step, $\hat{\pi}_0^E$ is computed using $\hat{\pi}_0^B$ or $\hat{\pi}_0^L$ as $\hat{\pi}_0^I$ and then in the next step, $\hat{\pi}_0^E$ is again computed using the previous value of $\hat{\pi}_0^E$ as $\hat{\pi}_0^I$. All the estimators discussed here may yield different numerical output for a given dataset but these estimators are only variants of the proposed estimation method.

\section{Expectation of $p$-values under the alternative}
For implementation of $\hat{\pi}_0^E$, appropriate estimate of $e$ is required. Probability density function $f_1^{\delta_i}(p)$ (for convenience this is written to be $f_{\delta_i}(p)$, henceforth) for each $p$-value, under the alternative with effect size $\delta_i$ is required to get explicit expression for $e_i=e_{\delta_i}$, $i\in\mathcal{F}$. The subscript $i$ present in the effect sizes are not specified in this section. For different testing scenarios, analytical expression of $e_{\delta}$ can be obtained using one of the following relations
\begin{equation}\label{eq:efQ}
    e_\delta=\int_0^1 p\,f_{\delta}(p)\, dp\quad \textrm{or} \quad e_\delta= \int_0^1 Q_\delta(p)\, dp. 
\end{equation}
\indent Let $X_1, X_2, ..., X_n$ be a random sample from a normal distribution with unknown mean $\mu$ and known variance $\sigma^2$. Consider the testing problem 
\begin{equation}\label{eq:hyp1}
    H_0: \mu=0\quad \quad \textrm{versus}\quad\quad H_1: \mu>0.
\end{equation}
Let $\bar{X}$ denote $(1/n)\sum_{i=1}^n X_i$, the sample mean. For the testing problem in \ref{eq:hyp1}, usual $Z$-test with right sided critical region is appropriate. The conventional test statistic $Z=\sqrt{n}\bar{X}/\sigma$ follows standard normal distribution, under $H_0$. The sample counterpart of $Z$ be $z$. The corresponding $p$-value is $p=P_{H_0}(Z>z)$. Under $H_1$, $Z$ is normally distributed with mean $\sqrt{n}\delta$ and variance 1. Here the effect size of the test $\delta=\mu/\sigma$. For the specified alternative hypothesis, $\delta$ cannot be negative. Thus, maximum likelihood estimator for $\delta$ is $\hat{\delta}=max\{0, \bar{X}/\sigma\}$. Obviously $\hat{\delta}$ is consistent for $\delta$, as required. One can obtain $\hat{e}_\delta$ by replacing $\delta$ by $\hat{\delta}$ in 
\begin{equation}\label{eq:test1edelta}
    e_{\delta}=1-E_{X\sim N(0,1)}\left\{\Phi(X+\sqrt{n}\delta)\right\}.
\end{equation}
Here, $\Phi$ denotes the cumulative distribution function of standard normal distribution and the expectation is taken with respect to $X$ from a standard normal distribution. For derivation of \ref{eq:test1edelta}, see Hung et al., 1997.\\

Now let $X_1, X_2, ..., X_n$ be a random sample from a normal distribution with unknown mean $\mu$ and unknown variance $\sigma^2$. Consider the testing problem 
\begin{equation}\label{eq:hyp2}
    H_0: \mu=0 \quad\quad \textrm{versus}\quad\quad H_1: \mu\neq 0.
\end{equation}
For the testing problem in \ref{eq:hyp2}, single sample two-sided $t$-test is appropriate with the test-statistic $T=\sqrt{n}\bar{X}/S$. Here $S^2$ is $(1/(n-1))\sum_{i=1}^n (X_i-\bar{X})^2$, the sample variance. $T$ follows the $t$-distribution with degrees of freedom (df) $n-1$, under $H_0$ while under $H_1$, it is distributed as the non-central $t$ with df $n-1$ and with non-centrality parameter (ncp) $\sqrt{n}\delta$. As in the earlier case, $\delta=\mu/\sigma$ and its maximum likelihood estimator is $\hat{\delta}=\sqrt{n}\bar{X}/S$. Let $F_{t_\nu}$ denote the cumulative distribution function of $t$-distribution with df $\nu$. Then the $p$-value corresponding to the two-sided test is defined as $p=2(1-F_{t_{n-1}}(|T|))$. Using the corresponding expression of $Q_\delta$, given from Cheng et al., 2015 in \ref{eq:efQ},
\begin{eqnarray*}
e_\delta &=&\int\limits_{0}^{1}[F_{t_{n-1,\sqrt{n}\delta}}(t_{n-1;\frac{p}{2}})- F_{t_{n-1,\sqrt{n}\delta}}(-t_{n-1;\frac{p}{2}})]\,dp\\
 &=&\int\limits_{0}^{1}F_{t_{n-1,\sqrt{n}\delta}}(t_{n-1;\frac{p}{2}})dp-\int\limits_{0}^{1}F_{t_{n-1,\sqrt{n}\delta}}(-t_{n-1;\frac{p}{2}})\,dp\\
 &=&I_1-I_2, \,\textrm{say}.\\
\end{eqnarray*}
Here, $F_{t_{\nu,\eta}}$ denotes cumulative distribution function $t$ distribution with df $\nu$ and ncp $\eta$ and $t_{\nu;x}$ denotes upper $x$-point of $t$ distribution with df $\nu$. Transforming $p$ to $v$ such that, $t_{n-1;p/2}=F_{t_{n-1;\sqrt{n}\delta}}^{-1}(1-p/2)=v$ we get 
\begin{equation*}
I_1= 2\int_0^\infty F_{t_{n-1,\sqrt{n}\delta}}(v)\, f_{t_{n-1,\sqrt{n}\delta}}(v)\,dv.
\end{equation*}
Here, $f_{t_{\nu,\eta}}$ denote density function of $t$ distribution with df $\nu$ and ncp $\eta$. Also let $t_{\nu}(a,b)$ denote truncated $t$ distribution with region of truncation $(a,b)$. Note that $t$ distribution is symmetric about $0$. Thus, the density function of $t_{\nu}(-\infty,0)$ is $(1/2)f_{t_{n-1}}$ for any df $\nu$. Therefore, 
\begin{eqnarray*}
I_1&=&\int\limits_{0}^{\infty}F_{t_{n-1,\sqrt{n}\delta}}(v)\frac{f_{t_{n-1}}(v)}{\frac{1}{2}}dv\\
&=&E_{X\sim t_{n-1}(-\infty ,0)}\{F_{t_{n-1,\sqrt{n}\delta}}(X)\}.
\end{eqnarray*} 
Here the expectation is taken with respect to $X$ from the specified truncated $t$ distribution. $I_2$ can be evaluated similarly. Hence from $e_\delta=I_1-I_2$,
\begin{equation}\label{eq:test2edelta}
    e_\delta=E_{X\sim t_{n-1}(-\infty,0)}\{F_{t_{n-1,\sqrt{n}\delta}}(X)\}-E_{X\sim t_{n-1}(0,\infty)}\{F_{t_{n-1,\sqrt{n}\delta}}(X)\}.
\end{equation}
\indent Consider $X_1, X_2, ..., X_{n_1}$ and $Y_1, Y_2, ..., Y_{n_2}$ be two random samples of size $n_1$ and $n_2$ respectively, from normal distribution with unknown mean $\mu_1$ and from normal distribution with unknown mean $\mu_2$. Assume that, the two normal populations have common variance $\sigma^2$. Consider the testing problem
\begin{equation}\label{eq:hyp3}
    H_0:\mu_1=\mu_2\quad\quad\textrm{versus}\quad\quad H_1:\mu_1\neq\mu_2.
\end{equation}
For the testing problem in \ref{eq:hyp3}, two sample two-sided  $t$-test is appropriate with $T=(\bar{X}-\bar{Y})/\sqrt{S^2(1/n_1+1/n_2)}$. Here, $S^2$ is $((n_1-1)S_X^2+(n_2-1)S_Y^2)/(n_1+n_2-2)$, the pooled sample variance ($S_X^2$ and $S_Y^2$ are sample variances of $X$-sample and $Y$-sample, respectively). In this case, $T$ follows $t$ distribution with df $n_1+n_2-2$ and ncp $\delta\sqrt{n^{*}}$, where $n^{*}=n_1n_2/(n_1+n_2)$. Here, $p=2(F_{t_{n_1+n_2-1}}(|T|))$. Using the corresponding $Q_\delta$ from Cheng et al., 2015 in \ref{eq:efQ} 
\begin{equation}\label{eq:etest3}
    e_\delta=E_{X\sim t_{n_1+n_2-2}(-\infty ,0)}\{F_{t_{n_1+n_2-2,\sqrt{n^*}\delta}}(X)\}-E_{X\sim t_{n_1+n_2-2}(0,\infty)}\{F_{t_{n_1+n_2-2,\sqrt{n^*}\delta}}(X)\}.
\end{equation}
Derivation of this result is similar to that of the result in \ref{eq:test2edelta}.

\section{Simulation study}
In order to generate an artificial dataset suitable for performing multiple single sample two-sided t-tests, the following steps have been carried out. Cheng et al.(2015) have established that, $\hat{\pi}_0^U$ performs better than some well-known estimators for $\pi_0$. Moreover, the set-up behind formulation of $\hat{\pi}_0^U$ is quite similar to that of the current work. Thus, $\hat{\pi}_0^U$ is a reasonable competitor of the proposed method of estimation. As mentioned in subsection 2.2, several variants of the proposed estimator $\hat{\pi}_0^E$ may be constructed and implemented to real data analysis. An effort has been made here to study relative performance of the following variants of $\hat{\pi}_0^E$ for simulated data.
\begin{eqnarray*}
&\hat{\pi}_0^{E_1}: \textrm{The proposed estimator with}\, \hat{\pi}_0^I=\hat{\pi}_0^B.\\
&\hat{\pi}_0^{E_2}: \textrm{The proposed estimator with}\, \hat{\pi}_0^I=\hat{\pi}_0^L.\\
&\hat{\pi}_0^{E_3}: \textrm{The proposed estimator with}\, \hat{\pi}_0^I=\hat{\pi}_0^{E_1}.\\
&\hat{\pi}_0^{E_4}: \textrm{The proposed estimator with}\, \hat{\pi}_0^I=\hat{\pi}_0^{E_2}.
\end{eqnarray*}
This may throw some light on the effect of choosing appropriate initial estimator.

\subsection{Simulation setting}
In order to generate an artificial dataset suitable for performing multiple single sample two-sided $t$-tests, the following steps have benn carried out. Set $m$, number of tests to be performed and $n$, size of the available sample observations for performing each test. For fixed $\pi_0= 0.1, 0.2, ..., 0.9$, find $m_0=[m\pi_0]$ and $m_1=m-m_0$. True mean under the null hypotheses, $\mu_0$ is set to be $0$. $\mathcal{T}$ is constructed by $m_0$ randomly chosen indices from $\mathcal{I}$. For the randomly generated $\mathcal{T}$, fix $\mu_i=\mu_0$ for $i\in \mathcal{T}$. Half of the $\mu$-values under alternative are generated from uniform$(0,0.5)$ distribution and the other half from uniform$(-0.5,0)$ distribution. For introducing weak correlation structure, two positive integers $b$ and $r$ are so chosen that, $m=b\times r$. Consider the covariance matrix
\begin{center}
   $\Sigma^{m\times m}=\left(\begin{array}{cccc}
\sigma_1^2\Sigma_\rho &0  &...  &0 \\ 
 0&\sigma_2^2\Sigma_\rho  & ... & 0\\ 
 \vdots&\vdots  &\vdots  &\vdots \\ 
 0&0  & ... & \sigma_r^2\Sigma_\rho
\end{array} \right)$. 
\end{center}
Here, $\sigma_i^2$'s are independently generated from exponential$(10)/3$ distribution for $i=1, 2, ..., r$. Explicit structure of $\Sigma_\rho$ is $\Sigma_\rho=(\rho^{|i-j|})$ for $i, j= 1, 2, ..., b$. For the $i$-th array of the dataset, a sample of size $n$ is generated from normal$(\mu_i,\Sigma_{i,i})$ for $i\in \mathcal{I}$. The numerical values for the flexible parameters in the given simulation studies are taken as $m=1000$, $n=25, 50$, $b=100$, $r=10$ and $\rho=0, 0.2, 0.5$.

\subsection{Simulation results}
The results reported in this subsection are based on $N=1000$ experiments for each of the simulation settings mentioned in the previous subsection. As discussed in subsection 2.1, Figure 1 is presented here to assess the strength of the bias correction. From figure 1, it is evident that, the behaviour of $e$ and $\tilde{e}$ are quite similar as functions of $\pi_0$. The difference between $e$ and $\tilde{e}$ reduces as $n$ increases. For each of the initial estimators, the performance of $\hat{e}$ gets better with increase in $n$. Performance of $\hat{e}$ is worst around $0.5$ and it improves for lower and higher values of $\pi_0$. It is best at $\pi_0=0.9$. It is also to be noted that, $\hat{e}$ approaches $e$ most effectively when the choice of initial estimator estimator is $\hat{\pi}_0^{B}$. This efficient way of reducing bias, at least for the larger values of $\pi_0$, makes the proposed method of estimation different from the conservative estimators for $\pi_0$. Now, the different variants of $\hat{\pi}_0^E$ and $\hat{\pi}_0^U$ are compared with respect to 
\begin{equation*}
    Bias(\hat{\pi}_0)=\frac{1}{N}\sum_{i=1}^N (\hat{\pi}_0^{(i)}-\pi_0)\quad\quad\textrm{and}\quad\quad MSE(\hat{\pi}_0)=\frac{1}{N}\sum_{i=1}^N (\hat{\pi}_0^{(i)}-\pi_0)^2.
\end{equation*}       
Here, $\hat{\pi}_0$ denotes any candidate estimator for $\pi_0$ and $\hat{\pi}_0^{(i)}$ is the observed value of that estimator in the $i$-th experiment, $i=1, 2, ..., N$. From Figure 2, it is evident that, the bias of each of the estimators considered in the simulation study decreases with increase in $n$ and $\pi_0$. For $\pi_0$ near $0.9$, bias of $\hat{\pi}_0^{E_1}$ and $\hat{\pi}_0^{U}$ is marginally negative. The behaviour of these two estimators are quite similar over the entire parameter space. It is to be noted that, though $\hat{\pi}_0^{E_1}$ has least bias for estimating $\pi_0$ (as seen from Figure 1), $\hat{\pi}_0^{E_2}$ seems to be safer for use in practice, as its bias is conservative even for higher values of $\pi_0$. Dependence among the rows of the data matrix seems to have no effect on the amount of bias. But, this is not the case for the mean squared error. From Figure 3, it is clear that, mean squared error suddenly increases for $\pi_0>0.5$ for $\rho=0.2$ and $0.5$. Mean squared errors of all the estimators under study, decreases with increasing $n$. As in case of bias, $\hat{\pi}_0^{E_1}$ and $\hat{\pi}_0^{U}$ behave more or less similarly, but $\hat{\pi}_0^{E_1}$ performs better than $\hat{\pi}_0^{U}$ for higher values of $\pi_0$. The opposite is true in case of smaller values of $\pi_0$. For $\pi_0>0.6$, $\hat{\pi}_0^{E_4}$ performs best under most of the simulation settings, but performance is relatively poor in the remaining portion of the parameter space. Even in this case, the performance of $\hat{\pi}_0^{E_2}$ is satisfactory and stable. Figure 5 is presented to get some insights of the distributions of the estimators. From the plot of empirical density functions, it is clear that the distribution of the estimators are skewed. Spread of the distributions of the estimators increase with $\rho$ and the mode of the estimators come closer to the true $\pi_0$ with increasing $\pi_0$. Comparative study on kurtosis of the estimators are given in Figure 4. 

\begin{figure}[H]
\begin{center}
\includegraphics[height=8in,width=6.7in,angle=0]{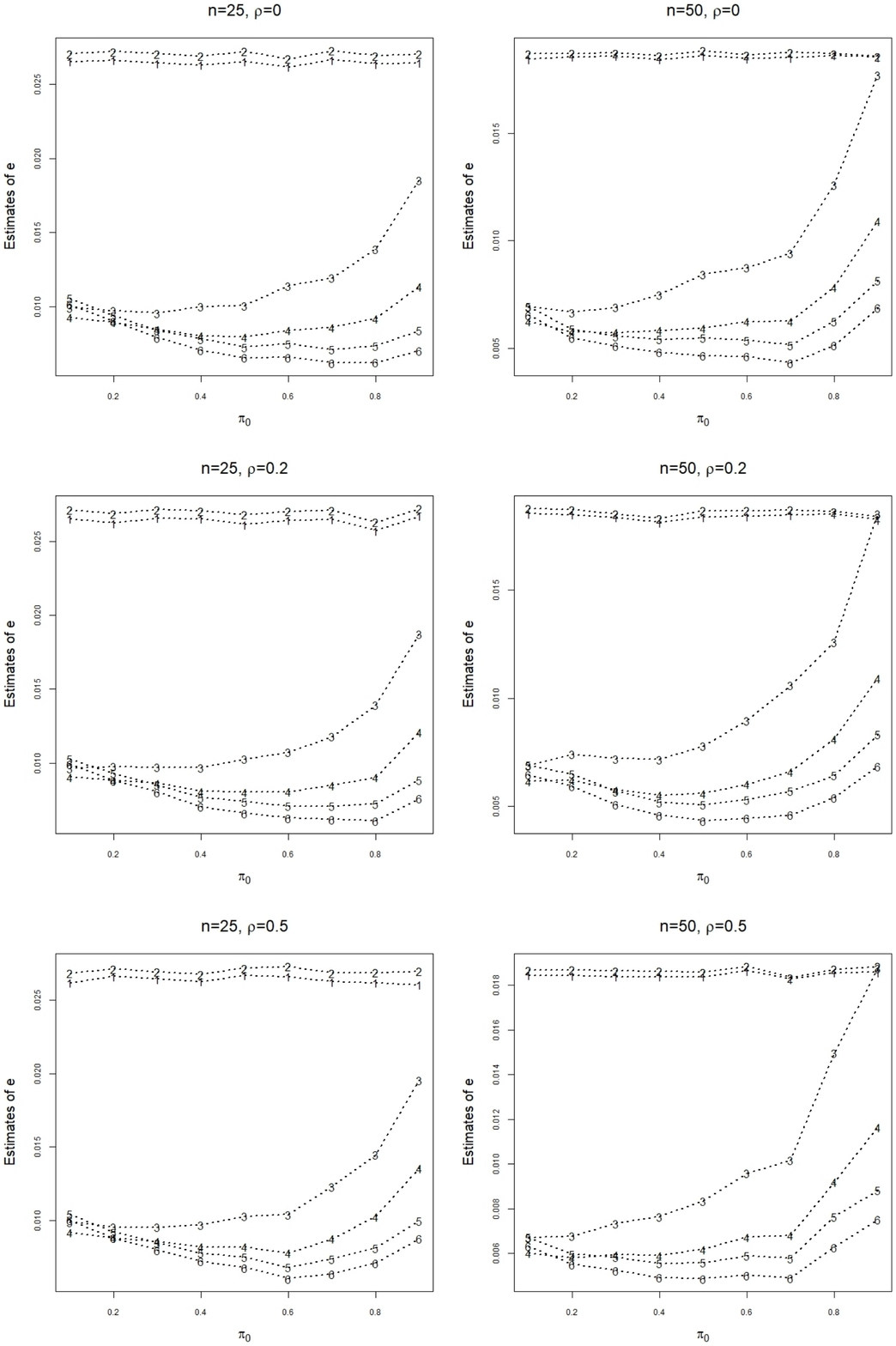}
\caption{Plots of $e$ (line 1), $\tilde{e}$ (line 2) and $\hat{e}$ using the initial estimators: $\hat{\pi}_0^{E_1}$ (line 3), $\hat{\pi}_0^{E_2}$ (line 4), $\hat{\pi}_0^{E_3}$ (line 5) and $\hat{\pi}_0^{4}$ (line 6) as functions of $\pi_0$.} 
\end{center}
\end{figure}

\begin{figure}[H]
\begin{center}
\includegraphics[height=9in,width=6.7in,angle=0]{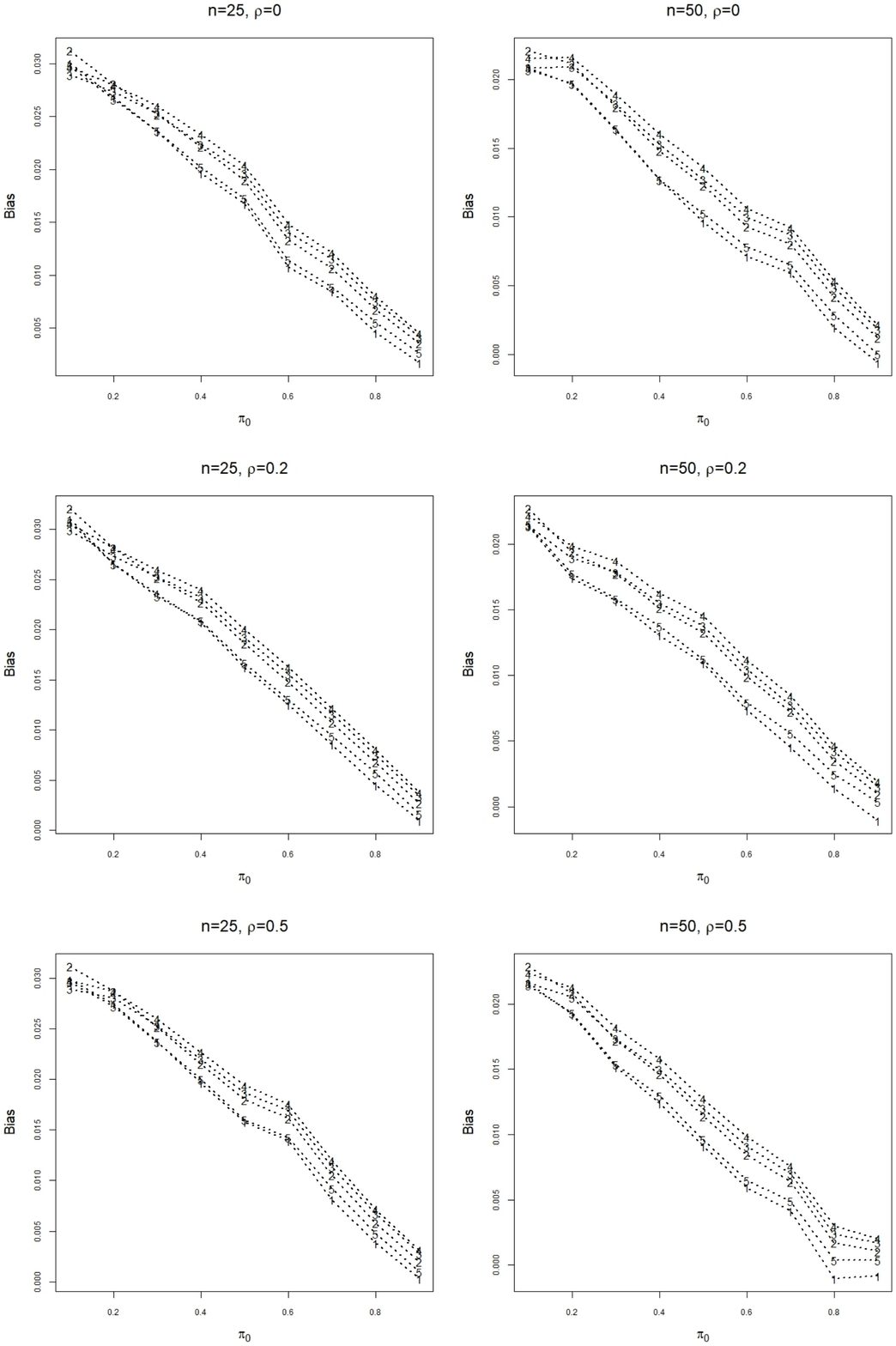}
\caption{Plots of bias for $\hat{\pi}_0^{E_1}$ (line 1), $\hat{\pi}_0^{E_2}$ (line 2), $\hat{\pi}_0^{E_3}$ (line 3), $\hat{\pi}_0^{E_4}$ (line 4) and $\hat{\pi}_0^{U}$ (line 5) as functions of $\pi_0$.} 
\end{center}
\end{figure}

\begin{figure}[H]
\begin{center}
\includegraphics[height=9in,width=6.7in,angle=0]{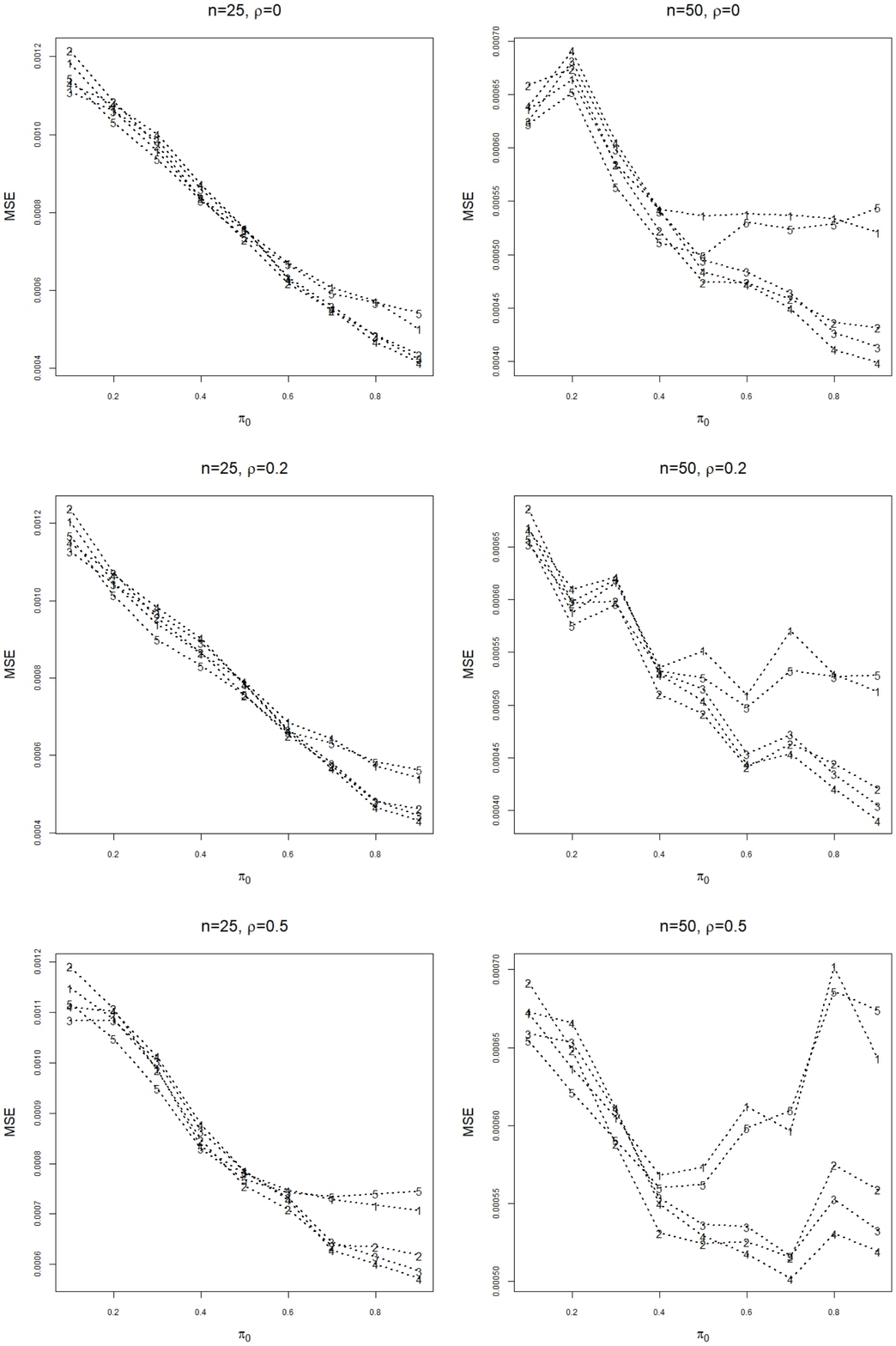}
\caption{Plots of mean squared error for $\hat{\pi}_0^{E_1}$ (line 1), $\hat{\pi}_0^{E_2}$ (line 2), $\hat{\pi}_0^{E_3}$ (line 3), $\hat{\pi}_0^{E_4}$ (line 4) and $\hat{\pi}_0^{U}$ (line 5) as functions of $\pi_0$.} 
\end{center}
\end{figure}

\begin{figure}[H]
\begin{center}
\includegraphics[height=9in,width=6.7in,angle=0]{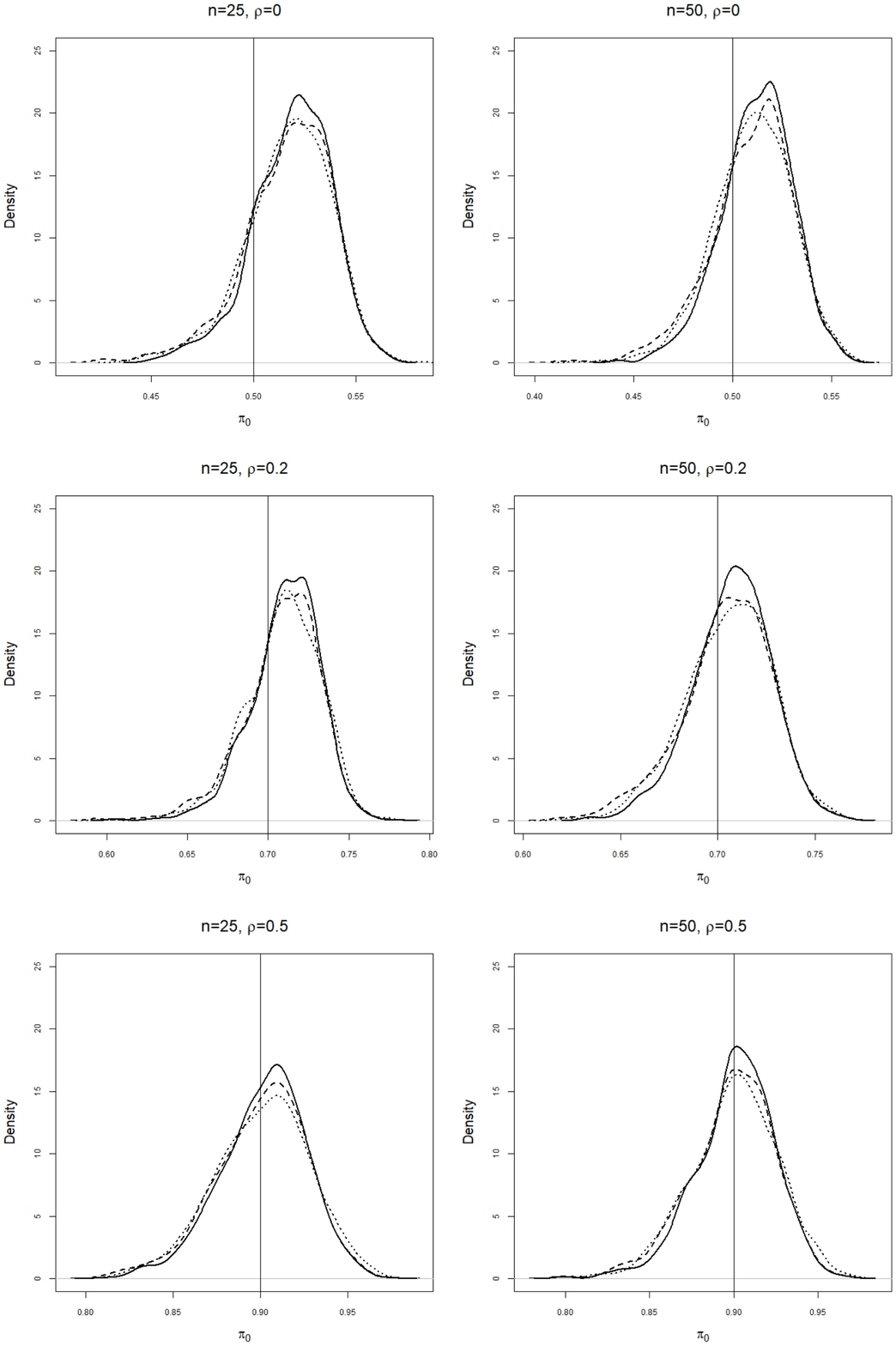}
\caption{Density plots of $\hat{\pi}_0^{E_1}$ (solid line), $\hat{\pi}_0^{E_2}$ (dashed line) and $\hat{\pi}_0^{U}$ (dotted line). True value of $\pi_0$ is indicated by the vertical line in each of the plots.} 
\end{center}
\end{figure}

\section{Data analysis}
For the purpose of case study, two popular datasets are used. The first one is the Leukemia dataset (see Golub et al., 1999) and the second is Prostate Cancer dataset (see Efron, 2012). In the first dataset, bone marrow samples are taken from 47 acute lymbhoblastic leukemia (ALL) patients and 25 acute myeloid leukemia (AML) patients. The samples were analysed using affymetrix arrays. There are 7128 genes in total. Objective of analysing this dataset is to estimate the proportion of genes which are significantly different among the two groups of patients: ALL and AML. In the second dataset, genetic expression levels for 6033 genes are obtained, for 50 normal control subjects and for 52 prostate cancer patients (see Efron, 2012). The objective of analysing this dataset is to estimate the proportion of differentially expressed genes. For both the datasets, two sample two-sided  t-tests are applicable. Different estimates of $\pi_0$ using the model based estimators are reported in Table 1.
\begin{table}[h!]
\centering
\caption{Different Estimates for the two datasets.}
\label{my-label}
\begin{tabular}{ccc}
ESTIMATES & LEUKEMIA DATA & PROSTATE DATA \\
$\hat{\pi}_0^{E_1}$ & 0.65192 & 0.90492 \\
$\hat{\pi}_0^{E_2}$  &0.65192   &0.90205\\
$\hat{\pi}_0^{E_3}$ &0.65192  &0.90457 \\ 
$\hat{\pi}_0^{E_4}$ &0.65192  &0.90463 \\
$\hat{\pi}_0^{U}$   & 0.62387 & 0.91258\\
\end{tabular}
\end{table}
It is to be noted that, all the variants of $\hat{\pi}_0^E$ yield the same estimate of $\pi_0$ for Leukemia data. This fact may be attributed to insensitivity of $\hat{\pi}_0^E$ to the choice of initial estimator. That is, initial estimates are so close to each other that, for all the cases $d$ turns out to be same. However, further investigation differs from the above explanation. In fact, the initial guesses for $m_1$ to be same is near impossible when we are testing a large number of hypotheses, however similar the $\hat{\pi}_0^I$'s may be. The reason behind having same estimates is that, there are a large number of false null hypotheses with strong signals. Scrutiny of intermediate steps while computing the estimates clearly shows that, the number of $\hat{e}_i$'s equals to $0$ is much larger than $d$, for any choice of $\hat{\pi}_0^I$. Similar statement holds good for $\hat{Q}_{\delta_i}(\lambda)$'s, for all $\lambda\in \Lambda$. This implies that, model based bias correction through $\hat{\pi}_0^E$ or $\hat{\pi}_0^U$ is not effective for datasets with a large number of strong signals. This being said, one should realize that, the task of identifying differentially expressed genes gets difficult when the signals are weak and in these kind situations, the proposed method of estimation may find its application. The Prostate data example is one such case.    
\section{Concluding remarks}
In the current work, a new method of estimating $\pi_0$ has been introduced and studied extensively when a large number of two-sided $t$-tests are performed, simultaneously. This scenario is very common in studies related to gene expression levels. The method of estimation is simple yet effective. The existing robust estimators for $\pi_0$ only utilize $p$-values from the corresponding multiple tests, whereas the proposed method of estimation also demands the original data. Thus, the proposed estimation procedure is applicable for datasets, which can be analyzed by simple methods like two-sample two-sided $t$-test. From the given simulation study, we identify the estimator proposed in Langaas et al., 2005 to be most effective as an initial estimator. However, this choice is not universal and further research on this issue is warranted. It is worth mentioning that, the one-step iterated estimators can also be used, but not recommended due to their huge computation time. The proposed estimator beats the estimator proposed in Cheng et al., 2015 over an important portion of the parameter space and reamins a viable alternative in the remaining part of the parameter space, even with the initial estimator being the estimator proposed in Storey, 2002. Thus, the proposed method of estimation may improve the existing literature. This work only concentrates on the estimation of the proportion of true null hypotheses. Exploring its performance for estimating the false discovery rate and construction of related adaptive algorithms for controlling the same remain an interesting study to be taken up in future.    

\section*{Reference}

\begin{footnotesize}
Benjamini, Y., \& Hochberg, Y. (1995). Controlling the false discovery rate: a practical and powerful approach to multiple testing. \textit{Journal of the royal statistical society. Series B (Methodological)}, 289-300.\\
\\
Cheng, Y., Gao, D., \& Tong, T. (2015). Bias and variance reduction in estimating the proportion of true-null hypotheses. \textit{Biostatistics}, 16(1), 189-204.\\
\\
Efron, B. (2012). Large-scale inference: \textit{empirical Bayes methods for estimation, testing, and prediction} (Vol. 1). Cambridge University Press.\\
\\
Finner, H., \& Gontscharuk, V. (2009). Controlling the familywise error rate with plug‐in estimator for the proportion of true null hypotheses. \textit{Journal of the Royal Statistical Society: Series B (Statistical Methodology)}, 71(5), 1031-1048.\\
\\
Golub, T. R., Slonim, D. K., Tamayo, P., Huard, C., Gaasenbeek, M., Mesirov, J. P., Coller, H., Loh, M.L., Downing, J.R., Caligiuri, M.A., Bloomfield, C. D., \& Lander, E.S. (1999). Molecular classification of cancer: class discovery and class prediction by gene expression monitoring. \textit{science}, 286(5439), 531-537.\\
\\
Hochberg, Y., \& Benjamini, Y. (1990). More powerful procedures for multiple significance testing. \textit{Statistics in medicine}, 9(7), 811-818.\\
\\
Hung, H. J., O'Neill, R. T., Bauer, P., \& Kohne, K. (1997). The behavior of the p-value when the alternative hypothesis is true. \textit{Biometrics}, 11-22.\\
\\
Jiang, H., \& Doerge, R. W. (2008). Estimating the proportion of true null hypotheses for multiple comparisons. \textit{Cancer informatics}, 6, 117693510800600001.\\
\\
Langaas, M., Lindqvist, B. H., \& Ferkingstad, E. (2005). Estimating the proportion of true null hypotheses, with application to DNA microarray data.\textit{ Journal of the Royal Statistical Society: Series B (Statistical Methodology)}, 67(4), 555-572.\\
\\
Miller, C. J., Genovese, C., Nichol, R. C., Wasserman, L., Connolly, A., Reichart, D., Hopkins, A., Schneider, A. \& Moore, A. (2001). Controlling the false-discovery rate in astrophysical data analysis. \textit{The Astronomical Journal}, 122(6), 3492.\\
\\
Qu, L., Nettleton, D., \& Dekkers, J. C. (2012). Improved Estimation of the Noncentrality Parameter Distribution from a Large Number of t‐Statistics, with Applications to False Discovery Rate Estimation in Microarray Data Analysis. \textit{Biometrics}, 68(4), 1178-1187.\\
\\
Storey, J. D. (2002). A direct approach to false discovery rates. \textit{Journal of the Royal Statistical Society: Series B (Statistical Methodology)}, 64(3), 479-498.\\
\\
Storey, J. D., Taylor, J. E., \& Siegmund, D. (2004). Strong control, conservative point estimation and simultaneous conservative consistency of false discovery rates: a unified approach. \textit{Journal of the Royal Statistical Society: Series B (Statistical Methodology)}, 66(1), 187-205.\\
\\
Storey, J. D., \& Tibshirani, R. (2003). SAM thresholding and false discovery rates for detecting differential gene expression in DNA microarrays. In \textit{The analysis of gene expression data} (pp. 272-290). Springer, New York, NY.\\
\\
Turkheimer, F. E., Smith, C. B., \& Schmidt, K. (2001). Estimation of the number of “true” null hypotheses in multivariate analysis of neuroimaging data. \textit{NeuroImage}, 13(5), 920-930.\\
\\
Wang, H. Q., Tuominen, L. K., \& Tsai, C. J. (2010). SLIM: a sliding linear model for estimating the proportion of true null hypotheses in datasets with dependence structures.\textit{ Bioinformatics}, 27(2), 225-231.

\end{footnotesize}
\end{document}